\documentclass{amsproc}

\newtheorem{theorem}{Theorem}[section]
\newtheorem{lemma}[theorem]{Lemma}

\theoremstyle{definition}
\newtheorem{definition}[theorem]{Definition}
\newtheorem{example}[theorem]{Example}
\newtheorem{problem}[theorem]{Problem}
\newtheorem{corollary}[theorem]{Corollary}
\newtheorem{proposition}[theorem]{Proposition}
\newtheorem{conjecture}[theorem]{Conjecture}

\theoremstyle{remark}

\numberwithin{equation}{section}

\newcommand{\Ff}{{\mathbb F}}
\newcommand{\Zz}{{\mathbb Z}}

\newcommand{\Qq}{{\mathbb Q}}

\newcommand{\Ll}{{\mathcal L}}

\newcommand{\cP}{{\mathcal P}}
\newcommand{\cB}{{\mathcal B}}
\newcommand{\cD}{{\mathcal D}}

\newcommand\MM{{\overline M}}   %

\def\Tr{\operatorname{Tr}}

\def\PG{\operatorname{PG}}
\def\AG{\operatorname{AG}}

\begin{document}

\title[$p$-ranks and Smith normal forms of $2$-designs]
{Recent results on $p$-ranks and Smith normal forms of some $2-(v,k,\lambda)$ designs}

\author{Qing Xiang}
\address{Department of Mathematical Sciences, University of
Delaware, Newark, DE 19716} \email{xiang@math.udel.edu}
\thanks{The author was supported in part by NSA grant MDA904-01-1-0036.}
\date{January 21, 2004}


\begin{abstract}
We survey recent results on $p$-ranks and Smith normal forms of
some $2-(v,k,\lambda)$ designs. In particular, we give a
description of the recent work in \cite{csx} on the Smith normal
forms of the 2-designs arising from projective and affine spaces
over $\Ff_q$.
\end{abstract}

\maketitle

\section{Introduction}
A {\it $2-(v,k,\lambda)$ design} is a pair $(\cP, {\mathcal B})$
where ${\mathcal P}$ is a $v$-set and ${\mathcal B}$ is a
collection of $b$ subsets of ${\mathcal P}$ ({\it blocks}), each
of size $k$, such that any 2-subset of ${\mathcal P}$ is contained
in exactly $\lambda$ blocks. Simple counting arguments show that
$b=\frac {\lambda v(v-1)} {k(k-1)},$ and the number of blocks
containing each element of $\cP$ is $\frac{\lambda (v-1)}{k-1}$,
which will be denoted by $r$ (called the {\it replication number}
of the design). The {\it order} of the 2-design, denoted by $n$,
is defined to be $r-\lambda$. A $2-(v,k,\lambda)$ design $(\cP,
{\mathcal B})$ is said to be {\it simple} if it does not have
repeated blocks (i.e., ${\mathcal B}$ is a set). We will only
consider simple 2-designs in this paper. A simple
$2-(v,k,\lambda)$ design $(\cP, {\mathcal B})$ is called {\it
symmetric}  (or {\it square}) if $b=v$.

Classical examples of 2-designs arise from projective and affine
spaces over finite fields. Let ${\rm PG}(m,q)$ be the
$m$-dimensional projective space over the finite field $\Ff_q$,
where $q$ is a prime power, let $\AG(m,q)$ be the $m$-dimensional
affine space over $\Ff_q$, and let $\begin{bmatrix} m\\
i\end{bmatrix}_q$ denote the number of $i$-dimensional subspaces
of an $m$-dimensional vector space over $\Ff_q$. We have the
following classical examples of 2-designs.

\begin{example}\label{projexample}
Let $m\geq 2$ and $m\geq d\geq 2$ be integers. The points of $\PG(m,q)$ and the $(d-1)$-dimensional subspaces of
$\PG(m,q)$ form a 2-design with parameters $v=\begin{bmatrix}m+1
\\1\end{bmatrix}_q=(q^{m+1}-1)/(q-1)$, $k=\begin{bmatrix}d\\ 1\end{bmatrix}_q=(q^d-1)/(q-1)$,
$r=\begin{bmatrix}m\\
d-1\end{bmatrix}_q$, $\lambda=\begin{bmatrix}m-1\\
d-2\end{bmatrix}_q$, and $b=\begin{bmatrix} m+1 \\
d\end{bmatrix}_q$.
\end{example}

\begin{example}\label{affexample}
Let $m\geq 2$ and $m-1\geq d\geq 1$ be integers. The points of $\AG(m,q)$ and the $d$-flats of $\AG(m,q)$ form a
2-design with parameters $v=q^m$, $k=q^d$, $r=\begin{bmatrix}
m\\d\end{bmatrix}_q$, $\lambda=\begin{bmatrix}m-1 \\
d-1\end{bmatrix}_q$, and  $b=q^{m-d}\begin{bmatrix}m \\
d\end{bmatrix}_q$. Here the $d$-flats of $\AG(m,q)$ are the cosets
of $d$-dimensional subspaces of the underlying $m$-dimensional
vector space over $\Ff_q$.
\end{example}

In this paper, we will consider the 2-designs in the above
examples in detail. Other designs considered are difference sets
and unitals, which will be defined in the coming sections.

Given two $2-(v,k,\lambda)$ designs $\cD_1=(\cP_1,\cB_1)$ and
$\cD_2=(\cP_2,\cB_2)$, we say that $\cD_1$ and $\cD_2$ are {\it
isomorphic} if there exists a bijection $\phi: \cP_1\rightarrow
\cP_2$ such that $\phi(\cB_1)=\cB_2$ and for all $p\in \cP_1$ and
$B\in\cB_1$, $p\in B$ if and only if $\phi(p)\in\phi(B)$. An {\it
automorphism} of a 2-design is an isomorphism of the design with
itself. The set of all automorphisms of a 2-design forms a group,
{\it the (full) automorphism group} of the design. An {\it
automorphism group} of a 2-design is any subgroup of the full
automorphism group.

Isomorphism of designs can also be defined more algebraically by
using incidence matrices of designs, which we define now. Let
$\cD=(\cP, {\mathcal B})$ be a $2-(v,k,\lambda)$ design and label
the points as $p_1,p_2,\ldots ,p_v$ and the blocks as
$B_1,B_2,\ldots ,B_b$. An {\it incidence matrix} of $(\cP,
{\mathcal B})$ is the matrix $A=(a_{ij})$ whose rows are indexed
by the blocks $B_i$ and whose columns are indexed by the points
$p_j$, where the entry $a_{ij}$ is 1 if $p_j\in B_i$, and 0
otherwise. From the definition of 2-designs, we see that the
matrix $A$ satisfies
\begin{equation}\label{defdesign}
A^{\top}A=(r-\lambda)I+\lambda J,\; AJ=kJ,
\end{equation}
where $I$ is the $v\times v$ identity matrix, and $J$ is the
$v\times v$ matrix with all entries equal to 1. Note that the
first equation in (\ref{defdesign}) tells us that the rank of $A$
over $\Qq$ is $v$. Now let $\cD_1=(\cP_1,\cB_1)$ and
$\cD_2=(\cP_2,\cB_2)$ be two $2-(v,k,\lambda)$ designs, and let
$A_1$ and $A_2$ be incidence matrices of $\cD_1$ and $\cD_2$
respectively. Then $\cD_1$ and $\cD_2$ are isomorphic if and only
if there are permutation matrices $P$ and $Q$ such that
\begin{equation}\label{isom}
PA_1Q=A_2,
\end{equation}
that is, the matrices $A_1$ and $A_2$ are permutation equivalent.

Next we define codes, $p$-ranks, and Smith normal forms of
2-designs. Let $\cD$ be a $2-(v,k,\lambda)$ design with incidence
matrix $A$. The {\it $p$-rank} of $\cD$ is defined as the rank of
$A$ over a field $F$ of characteristic $p$, and it will be denoted
by ${\rm rank}_p(\cD)$. The $F$-vector space spanned by the rows
of $A$ is called the {\it (block) code} of $\cD$ over $F$, which
is denoted by $C_F(\cD)$. If $F=\Ff_q$, where $q$ is a power of
$p$, then we denote the code of $\cD$ over $\Ff_q$ by $C_q(\cD)$.
We proceed to define the Smith normal form of $\cD$. Let $R$ be a
principal ideal domain. Viewing $A$ as a matrix with entries in
$R$, we can find (see for example, \cite{cohn}) two invertible
matrices $U$ and $V$ over $R$ such that
\begin{equation}\label{snf}
UAV=\begin{pmatrix}
d_1 & 0 & 0 & \cdots & 0\\
0 & d_2 & 0\\
0 & & \ddots & & \vdots\\
\vdots & & & d_{v-1}& 0\\
0 & & \cdots & 0 & d_v\\
0 & & \cdots & & 0\\
\vdots & & \ddots& & \vdots\\
0 & & \cdots & & 0
\end{pmatrix}
\end{equation}
with $d_1|d_2|\cdots |d_v$. The $d_i$ are unique up to units in
$R$. When $R=\Zz$, the $d_i$ are integers, and they are called the
{\it invariant factors of $A$}; the matrix in the right hand side
of (\ref{snf}) (now with integer entries) is called the {\it Smith
normal form {\rm (SNF)} of $A$}. We define the {\it Smith normal
form of $\cD$} to be that of $A$. Smith normal forms and $p$-ranks
of 2-designs can help distinguish nonisomorphic 2-designs with the
same parameters: let $\cD_1=(\cP_1,\cB_1)$ and
$\cD_2=(\cP_2,\cB_2)$ be two $2-(v,k,\lambda)$ designs with
incidence matrices $A_1$ and $A_2$ respectively. From (\ref{isom})
we see that if $\cD_1$ and $\cD_2$ are isomorphic, then $A_1$ and
$A_2$ have the same Smith normal form over $\Zz$; hence $\cD_1$
and $\cD_2$ have the same Smith normal form, in particular, ${\rm
rank}_p(\cD_1)={\rm rank}_p(\cD_2)$ for any prime $p$. Furthermore
Smith normal forms of symmetric designs are used by Lander
\cite{lander} to construct a sequence of $p$-ary codes which carry
information about the designs. Therefore it is interesting to
study the codes and Smith normal forms of 2-designs.

We collect some general results on $p$-ranks and Smith normal
forms of 2-designs.

\begin{theorem}[Theorem 2.4.1 in \cite{assmus}]\label{genrank}
Let $\cD=(\cP, \cB)$ be a $2-(v,k,\lambda)$ design with order $n$.
Let $p$ be a prime and let $F$ be a field of characteristic $p$,
where $p\not {|}n$. Then
$${\rm rank}_p(\cD)\geq v-1$$
with equality if and only if $p|k$; in the case of equality we
have that $C_F(\cD)=(F{\bf j})^{\perp}$ and otherwise
$C_F(\cD)=F^{\cP}$. Here ${\bf j}$ is the all one row vector of
length $v$.
\end{theorem}

For primes $p$ dividing $n$, Klemm \cite{klemm1} proved the
following result.

\begin{theorem}[Klemm]\label{klemmrank}
Let $\cD$ be a $2-(v,k,\lambda)$ design with order $n$ and let $p$
be a prime dividing $n$. Then
\begin{equation}\label{uppbound}
{\rm rank}_p(\cD)\leq \frac {b+1} {2}.
\end{equation}
Moreover, if $p\not{|}\lambda$ and $p^2\not{|}n$, then
$$C_p(\cD)^{\perp}\subseteq C_p(\cD)$$
and ${\rm rank}_p(\cD)\geq v/2$.
\end{theorem}

If the design $\cD$ above is symmetric (i.e., $b=v$), then the
bound in (\ref{uppbound}) becomes
$${\rm rank}_p(\cD)\leq \frac {v+1} {2}.$$
This bound is best possible. For example, any skew Hadamard design
with parameters $(4n-1, 2n-1, n-1)$ has $p$-rank equal to $2n$,
where $p$ is any prime divisor of $n$ \cite{mich}. However, for a
2-design which is not symmetric, $b$ is usually much larger than
$v$, and the bound in (\ref{uppbound}) becomes very weak.

Klemm \cite{klemm2} also proved some general results on Smith
normal forms of 2-designs. He actually proved his results for what
he called semi-block designs, in which the block sizes may not be
uniform. Here we restrict our attention to 2-designs only.

\begin{theorem}[Klemm]\label{klemmsnf}
Let $\cD$ be a $2-(v,k,\lambda)$ design with order $n$ and let
$t=\gcd(n,\lambda)$. Let $d_1, d_2, \dots , d_v$ be the invariant
factors of $\cD$, where $d_1|d_2|\cdots |d_v$. Then
\begin{enumerate}
\item $d_1=1$, $(d_1d_2\cdots d_i)^2|tn^{i-1}$ for $2\leq i\leq
v-1$, and $$(d_1d_2\cdots d_v)^2|(n+\lambda v)n^{v-1}.$$

\item $d_v|(rn/t)$, and $d_i|n$ for $2\leq i\leq v-1$.

\item $p|d_i$ for $(b+1)/2 < i\leq v$ and every prime $p$
dividing $n$.
\end{enumerate}
\end{theorem}

For symmetric designs, the above theorem can be improved.

\begin{theorem}[Klemm]\label{klemmsym}
Let $\cD$ be a $2-(v,k,\lambda)$ symmetric design with order $n$
and let $t=\gcd(n,\lambda)$. Let $d_1, d_2, \dots , d_v$ be the
invariant factors of $\cD$, where $d_1|d_2|\cdots |d_v$. Then
\begin{enumerate}
\item $d_1d_2\cdots d_v=kn^{(v-1)/2}$.
\item $d_v=kn/t$.
\item Let $p$ be a prime dividing $n$ and $p\not{|}\lambda$.
For any integer $x$, let $x_p$ denote the $p$-part of $x$. Then
$(d_id_{v+2-i})_p=n_p$, for $3\leq i\leq v-1$.
\end{enumerate}
\end{theorem}

The rest of the paper is organized as follows. In Section 2, we
consider two families of recently constructed cyclic difference
sets with parameters $((3^{m}-1)/2, 3^{m-1}, 2\cdot 3^{m-2})$.
These two families of difference sets have the same 3-ranks. Yet
they are shown in \cite{cx} to be inequivalent by using the
numbers of 3's in their SNF. In Section 3, we report the recent
work in \cite{csx} on the SNF of the designs in
Example~\ref{projexample} and \ref{affexample}. In Section 4, we
collect some recent results on $p$-ranks and SNF of unitals. We
did not intend to be comprehensive. So there might be some recent
results on or related to $p$-ranks and SNF of 2-designs not
mentioned here.

\section{The invariant factors of some cyclic difference sets}

Let $\cD=(\cP, \cB)$ be a $2-(v,k,\lambda)$ symmetric design with
a sharply transitive automorphism group $G$. Then we can identify
the elements of $\cP$ with the elements of $G$. After this
identification, each block of $\cD$ is now a $k$-subset of $G$.
Since $G$ acts sharply transitively on $\cB$, we may choose a base
block $D\subset G$, all other blocks in $\cB$ are simply
``translates'' $gD=\{gx\mid x\in D\}$ of $D$, where $g\in G$ and
$g\neq 1$. That $\cD$ is a symmetric design implies
$$|D\cap gD|=\lambda,$$
for all nonidentity element $g\in G$. That is, every nonidentity
element $g\in G$ can be written as $xy^{-1}$, $x,y\in D$, in
$\lambda$ ways. This leads to the definition of difference sets.

\begin{definition}
Let $G$ be a finite (multiplicative) group of order $v$. A
$k$-element subset $D$ of $G$ is called a {\em
$(v,k,\lambda)$-difference set\/} in $G$ if the list of
``differences'' $xy^{-1}$, $x,y\in D$, $x\neq y$, represents each
nonidentity element in $G$ exactly $\lambda$ times.
\end{definition}

In the above, we see that sharply transitive symmetric designs
give rise to difference sets. In the other direction, if $D$ is a
$(v,k,\lambda)$-difference set in a group $G$, then we can use the
elements of $G$ as points, and use the ``translates" $gD$ of $D$,
$g\in G$, as blocks, and we obtain a symmetric design $(G,
\{gD\mid g\in G\})$ with a sharply transitive automorphism group
$G$. (This design is usually called {\it the symmetric design
developed from $D$}, and will be denoted by ${\rm Dev}(D)$.) Hence
difference sets and sharply transitive symmetric designs are the
same objects.

We will only consider difference sets in abelian groups. Let $D_1$
and $D_2$ be two $(v,k,\lambda)$-difference sets in an abelian
group $G$. We say that $D_1$ and $D_2$ are {\em equivalent\/} if
there exists an automorphism $\sigma$ of $G$ and an element $g\in
G$ such that $\sigma(D_1)=D_2g$. Note that if $D_1$ and $D_2$ are
equivalent, then ${\rm Dev}(D_1)$ and ${\rm Dev}(D_2)$ are
isomorphic. Therefore one way to distinguish inequivalent
difference sets is to show that the symmetric designs developed
from them are nonisomorphic. This will be the approach we take in
this paper. For this purpose, we define $p$-ranks, invariant
factors, and the Smith normal form of a $(v,k,\lambda)$-difference
set $D$ to be that of the associated design ${\rm Dev}(D)$.

In the study of abelian difference sets, characters play an
important role. The following is a basic lemma in this area, see
\cite{turyn}.

\begin{lemma}\label{charlemma}
Let $G$ be an abelian group of order $v$ and let $D$ be a
$k$-subset of $G$. Then $D$ is a $(v,k,\lambda)$-difference set in
$G$ if and only if
\begin{equation}
\chi(D)\overline{\chi(D)}= {k-\lambda}
\end{equation}
for every nontrivial complex character $\chi$ of $G$. Here
$\chi(D)$ stands for $\sum _{d\in D} ^{}\chi(d)$.
\end{lemma}

We will see later that if $\gcd(v,k-\lambda)=1$, then the
computations of $p$-ranks and invariant factors of a difference
set $D$ in an abelian group $G$ depend on our understanding of the
algebraic integers $\chi(D)$, where $\chi$ runs through the
character group of $G$.

The classical examples of difference sets are the Singer
difference sets. They arise from the case $d=m$ in
Example~\ref{projexample}. We state this formally below.

\begin{theorem}\label{sing}
Let $q$ be a prime power and $m>1$ an integer. Then the points of
${\rm PG}(m,q)$ and the $(m-1)$-dimensional subspaces
(hyperplanes) of ${\rm PG}(m,q)$ form a symmetric design admitting
a cyclic sharply transitive automorphism group. That is, the
point-hyperplane design in ${\rm PG}(m,q)$ is developed from a
cyclic difference sets with parameters
\begin{equation}\label{classic}
v=\frac {q^{m+1}-1} {q-1},\; k=\frac {q^{m}-1}
{q-1},\;\lambda=\frac {q^{m-1}-1} {q-1}.
\end{equation}
\end{theorem}

The $p$-ranks of the Singer difference sets were known from 1960's
(see \cite{xiang} for detailed references). The Smith normal forms
of the Singer difference sets were worked out in full generality
by Sin \cite{sin}, and independently by Liebler \cite{liebler}.
Since the result on the SNF of Singer difference sets is a special
case of the more general results in Section 3, we will not state
their results here.

The parameters in (\ref{classic}) or the complementary parameters
of (\ref{classic}) are called {\it classical parameters}. It is
known that there are many infinite families of cyclic difference
sets with classical parameters which are inequivalent to the
Singer difference sets. For a survey of results up to 1999. we
refer the reader to \cite{xiang}. Most of the examples of cyclic
difference sets with classical parameters in that survey
\cite{xiang} have even $q$, where $q$ is as in (\ref{classic}).
More examples with odd $q$ were discovered recently. Here are two
examples where $q$ is a power of $3$.

We will use standard notation: $\Ff_{q^m}$ denotes the finite
field with $q^m$ elements, $\Ff_{q^m}^*$ is the multiplicative
group of $\Ff_{q^m}$, $\Tr
 _{q^m/q}$ denotes the trace from $\Ff_{q^m}$ to $\Ff_q$, and the map
 $\rho: \Ff_{q^m}^*\rightarrow \Ff_{q^{m}}^*/\Ff_q^*$ denotes the natural
epimorphism.

\begin{example}
 Let $q=3^e$, $e\geq 1$, let $m=3k$, $k$ a positive integer,
$d=q^{2k}-q^k+1$, and set
\begin{equation}\label{defD}
R=\{x\in \Ff_{q^{m}}\mid \Tr_{q^{m}/q}(x+x^d)=1\}.
\end{equation}
 Then $\rho(R)$ is a $((q^{m}-1)/(q-1), q^{m-1}, q^{m-2}(q-1))$ difference
 set in $\Ff_{q^{m}}^*/\Ff_q^*$.
\end{example}
This is proved by using the language of sequences with ideal
2-level autocorrelation in \cite{hkm} in the case $q=3$. See
\cite{cx} for a complete proof of this fact (the paper \cite{cx}
also showed that $R$ is a relative difference set). For
convenience, we will call this difference set $\rho(R)$ the HKM
difference set.

\begin{example}
Let $m\ge 3$ be an odd integer, let $d=2\cdot 3^{(m-1)/2}+1$, and
set
\begin{equation}\label{defD1}
R=\{x\in \Ff_{3^{m}}\mid \Tr_{3^{m}/3}(x+x^d)=1\}.
\end{equation}
 Then $\rho(R)$ is a $((3^{m}-1)/2, 3^{m-1}, 2\cdot 3^{m-2})$ difference
 set in $\Ff_{3^{m}}^*/\Ff_3^*$.
\end{example}
This was conjectured by Lin, and recently proved by Arasu, Dillon
and Player \cite{a}. For convenience, we will call this difference
set $\rho(R)$ the Lin difference set.

In the case $q=3$, $m=3k$, $k>1$, the $3$-rank of the HKM
difference set is $2m^2-2m$. This was shown in \cite{cx} and
\cite{nsh}. One can similarly show that the Lin difference set has
3-rank $2m^2-2m$, where $m>3$ is odd, see \cite{nsh}. Therefore
when $m$ is an odd multiple of 3, these two difference sets have
the same 3-rank. Hence they can not be distinguished by 3-ranks.
It is therefore natural to consider using the SNF of these two
families of difference sets to distinguish them. We first state a
lemma which is very useful for determining the SNF of
$(v,k,\lambda)$-difference sets with $\gcd(v,n)=1$.

\begin{lemma}\label{smithgloblem}
 Let $G$ be an abelian group of order $v$, let $p$ be a prime not dividing
 $v$, and let $\mathfrak{P}$ be a prime ideal in $\Zz[\xi_{v}]$ lying
 above $p$ , where $\xi_v$ is a complex primitive $v^{\rm th}$ root of unity. Let
 $D$ be a $(v,k,\lambda)$ difference set in $G$, and let $\alpha$ be a
 positive integer. Then the number of invariant factors of $D$ which are
 not divisible by $p^{\alpha}$ is equal to the number of complex
 characters $\chi$ of $G$ such that $\chi(D)\not\equiv 0 \; ({\rm mod} \;
\mathfrak{P}^{\alpha})$.
\end{lemma}

Setting $\alpha=1$ in Lemma~\ref{smithgloblem}, we see that the
$p$-rank of $D$ is equal to the number of complex characters
$\chi$ such that $\chi(D)\not\equiv 0$ (mod $\mathfrak{P}$). This
was proved by MacWilliams and Mann \cite{mm}. For a full proof of
the lemma, see \cite{cx}.

Using Lemma~\ref{smithgloblem}, Fourier transforms, and
Stickelberger's congruence on Gauss sums, we \cite{cx} computed
the number of 3's in the SNF of the Lin and HKM difference sets.

\begin{theorem}

Let $m>9$. Then the number of 3's in the Smith normal form of the
HKM difference sets with parameters $((3^{m}-1)/2, 3^{m-1}, 2\cdot
3^{m-2})$ is
$$\frac23m^4-4m^3-\frac{28}3m^2+62m+\epsilon(m)\cdot m.$$

The number of 3's in the Smith normal form of the Lin difference
sets when $m>7$ is
$$\frac23m^4-4m^3-\frac{14}3m^2+39m
+\delta (m) \cdot m.$$

The values of $\epsilon(m)$ and $\delta(m)$ are 0 or 1.
\end{theorem}

Based on numerical evidence, we conjecture that $\delta$ and
$\epsilon$ above are always 1. By direct calculations (i.e., not
using Gauss sums), the Smith normal form of the Lin difference set
with $m=9$ is:
$$1^{144}3^{1440}9^{1572}27^{1764}81^{1764}243^{1572}729^{1440}
2187^{144}6561^{1},$$
 where for example, $3^{1440}$ means the number of invariant factors of the Lin
difference set which are 3 is $1440$. The Smith normal form of the
HKM difference set with $m=9$ is:
$$1^{144}3^{1251}9^{1842}27^{1683}81^{1683}243^{1842}729^{1251}2187^{144}6561^1.$$
These computations were done by Saunders \cite{saunders}.

Since the two ``almost'' polynomial functions in Theorem 2.7 are
never equal when $m>9$, and since the Smith normal forms of the
Lin and HKM difference sets are also different when $m=9$, we have
the following conclusion:

\begin{theorem} Let $m$ be an odd multiple of 3. The Lin and HKM difference
sets with parameters
 $(\frac {3^m-1} {2}, 3^{m-1}, 2\cdot 3^{m-2})$ are inequivalent when $m>3$,
and the associated symmetric designs are nonisomorphic when $m>3$.
\end{theorem}

Therefore we successfully distinguished the HKM and Lin difference
sets by using the number of 3's in their SNF. At this point, it is
natural to ask whether it is true that two symmetric designs with
the same parameters and having the same SNF are necessarily
isomorphic. The answer to this question is negative. It is known
\cite{as} that the Smith normal form of a projective plane of
order $p^2$, $p$ prime, is
$$1^{r}p^{(p^4+p^2-2r+2)}(p^2)^{(r-2)}((p^2+1)p^2)^1,$$
where the exponents indicate the multiplicities of the invariant
factors and $r$ is the $p$-rank of the plane. (This also follows
from Theorem~\ref{klemmsym}.) That is, the $p$-rank of the plane
completely determines the Smith normal form of the plane. There
are four projective planes of order 9. The desarguesian one has
3-rank 37, while the other three all have 3-rank $41$ (cf.
\cite{sachar}), so the three non-desarguesian projective planes
have the same Smith normal form, yet they are nonisomorphic.
However, the answer to the following more restricted question is
not known.

\begin{problem}
If two cyclic difference sets with classical parameters have the
same Smith normal form, are the associated designs necessarily
isomorphic?
\end{problem}

\section{The invariant factors of the incidence matrices of points
and subspaces in ${\rm PG}(m,q)$ and ${\rm AG}(m,q)$}

In this section, we describe the recent work in \cite{csx} on the
SNF of the designs in Example~\ref{projexample}
and~\ref{affexample}. We will concentrate on the design coming
from projective geometry first. The SNF of the design coming from
${\rm AG}(m,q)$ follows from the results in the projective case.

Let ${\rm PG}(m,q)$ be the $m$-dimensional projective space over
$\Ff_q$ and let $V$ be the underlying $(m+1)$-dimensional vector
space over $\Ff_q$, where $q=p^t$, $p$ is a prime. For any $d$,
$1\leq d\leq m$, we will refer to $d$-dimensional subspaces of $V$
as $d$-subspaces and denote the set of these subspaces in $V$ as
$\mathcal L_{d}$. The set of projective points is then $\mathcal
L_1$. The pair $(\mathcal L_1, \mathcal L_{d})$, where $d>1$, with
incidence being set inclusion, is the 2-design in
Example~\ref{projexample}. Let $A$ be an incidence matrix of the
2-design $(\mathcal L_1, \mathcal L_{d})$.
So $A$ is a $b\times v$ $(0,1)$-matrix, where $b=\begin{bmatrix} m+1 \\
d\end{bmatrix}_q$ and $v=\begin{bmatrix}m+1
\\1\end{bmatrix}_q$. We will determine the Smith normal form of $A$.
There is a somewhat long history of this problem. We refer the
reader to \cite{csx} for a detailed account.

The following theorem shows that all but one invariant factor of
$A$ are $p$ powers.

\begin{theorem}\label{p'part}
Let $A$ be the matrix defined as above. The invariant factors of
$A$ are all $p$-powers except for the $v^\mathrm{th}$ invariant,
which is a $p$-power times $(q^d-1)/(q-1)$.
\end{theorem}

This was known at least from \cite{sin1}. For a detailed proof,
see \cite{csx}. In view of Theorem~\ref{p'part}, to determine the
SNF of $A$, it suffices to determine the multiplicity of $p^i$
appearing as an invariant factor of $A$. It will be convenient to
view $A$ as a matrix with entries from a $p$-adic local ring $R$
(some extension ring of $\Zz_p$, the ring of $p$-adic integers).
We will define this ring $R$ and introduce a sequence of
$R$-modules and a sequence of $q$-ary codes in the following
subsection.

\subsection{$R$-modules and $q$-ary codes}

Let $q=p^t$ and let $K=\Qq_p(\xi_{q-1})$ be the unique unramified
extension of degree $t$ over $\Qq_p$, the field of $p$-adic
numbers, where $\xi_{q-1}$ is a primitive $(q-1)^{\rm th}$ root of
unity in $K$. Let $R=\Zz_p[\xi_{q-1}]$ be the ring of integers in
$K$ and let $\mathfrak p$ be the unique maximal ideal in $R$ (in
fact, $\mathfrak p=pR$). Then $R$ is a principal ideal domain, and
the reduction of $R \ (\mathrm{mod}\,\mathfrak p)$ is $\Ff_q$.
Define $\bar x$ to be $x\ (\mathrm{mod}\,\mathfrak p)$ for $x\in
R$.

We now view the above matrix $A$ as a matrix with entries from
$R$. Define
$$
M_i=\{x\in R^{\mathcal L_1}\ \mid Ax^{\top}\in p^iR^{\mathcal
L_d}\}, \quad i=0,1,...
$$
Here we are thinking of elements of $R^{\mathcal L_1}$ as row
vectors of length $v$. Then we have a sequence of nested
$R$-modules
$$
R^{\mathcal L_1}=M_0\supseteq M_1\supseteq \cdots
$$
Define $\MM_i=\{({\bar x_1},{\bar x_2},\ldots ,{\bar x_v})\in
\Ff_q^{\mathcal L_1}\mid (x_1,x_2,\ldots ,x_v)\in M_i\}$, for
$i=0,1,2,\ldots$. For example,
\begin{equation}
\MM_1=\{({\bar x_1},{\bar x_2},\ldots ,{\bar x_v})\in
\Ff_q^{\mathcal L_1}\mid A\begin{pmatrix}
x_1\\x_2\\\vdots\\x_v\end{pmatrix}\in pR^{\mathcal L_d}\}.
\end{equation}
That is, $\MM_1$ is the dual code of the $q$-ary (block) code of
the 2-design $({\mathcal L_1}, {\mathcal L_d})$. We have a
sequence of nested $q$-ary code
$$
\Ff_q^{\mathcal L_1}=\MM_0\supseteq \MM_1\supseteq \cdots
$$
This is similar to what Lander did for symmetric designs, see
\cite{lander} and \cite[p.~399]{vanlintwilson}. Note that if
$i>\nu_p(d_v)$, where $\nu_p$ is the $p$-adic valuation and $d_v$
is the $v^{\rm th}$ invariant factor of $A$, then $\MM_i=\{0\}$.
It follows that there exists a smallest index $\ell$ such that
$\MM_{\ell}=\{0\}$. So we have a finite filtration
$$
\Ff_q^{\mathcal L_1}=\MM_0\supseteq \MM_1\supseteq\cdots\supseteq \MM_{\ell}=\{0\}.
$$
We have the following easy but important lemma. See \cite{csx} for
its proof.

\begin{lemma}\label{multiplicity}
For $0\leq i\leq {\ell-1}$, $p^i$ is an invariant factor of $A$
with multiplicity $\dim_{\Ff_q} (\MM_i/\MM_{i+1})$.
\end{lemma}

In what follows, we will determine $\dim_{\Ff_q}(\MM_i)$, for each
$i\geq 0$. In fact, we will construct an $\Ff_q$-basis for each
$\MM_i$. To this end, we construct a basis of $\Ff_q^{\mathcal
L_1}$ first.

\subsection{Monomial basis of $\Ff_q^{\mathcal L_l}$ and types of basis monomials}

Let $V=\Ff_q^{m+1}$. Then $V$ has a standard basis $v_0,v_1,\ldots
,v_{m}$, where
$$v_i=(\underbrace{0,0,\ldots ,0,1}_{i+1},0,\ldots ,0).$$
We regard $\Ff_q^V$ as the space of functions from $V$ to $\Ff_q$.
Any function $f\in \Ff_q^V$ can be given as a polynomial function
of $m+1$ variables corresponding to the $m+1$ coordinate
positions: write the vector ${\mathbf x}\in V$ as
$${\mathbf x}=(x_0,x_1,\ldots ,x_m)=\sum_{i=0}^{m}x_iv_i;$$
then $f=f(x_0,x_1,\ldots ,x_m)$. The function $x_i$ is, for
example, the linear functional that projects a vector in $V$ onto
its $i^{\rm th}$ coordinate in the standard basis.

As a function on $V$, $x_i^q=x_i$, for each $i=0,1,\ldots ,m$, so
we obtain all the functions via the $q^{m+1}$ monomial functions
\begin{equation}\label{affbasis}
\{\prod_{i=0}^mx_i^{b_i}\mid 0\leq b_i<q, i=0,1,\ldots ,m\}.
\end{equation}
Since the characteristic function of $\{0\}$ in $V$ is
$\prod_{i=0}^m(1-x_i^{q-1})$, we obtain a basis for
$\Ff_q^{V\setminus\{0\}}$ by excluding $x_0^{q-1}x_1^{q-1}\cdots
x_m^{q-1}$ from the set in (\ref{affbasis}) (some authors prefer
to exclude $x_0^0x_1^0\cdots x_m^0$, see \cite{glyn}).

The functions on $V\setminus\{0\}$ which descend to $\mathcal L_1$
are exactly those which are invariant under scalar multiplication
by $\Ff_q^*$. Therefore we obtain a basis ${\mathcal M}$ of
$\Ff_q^{\mathcal L_1}$ as follows.
$${\mathcal M}=\{\prod_{i=0}^mx_i^{b_i}\mid 0\leq b_i<q,
\sum_i b_i\equiv 0\; ({\rm mod}\; q-1), (b_0,b_1,\ldots ,b_m)\neq
(q-1,q-1,\ldots ,q-1)\}.$$ This basis ${\mathcal M}$ will be
called the {\it monomial basis} of $\Ff_q^{\mathcal L_1}$, and its
elements are called {\it basis monomials}.

Next we define the type of a nonconstant basis monomial. Let
$\mathcal{H}$ denote the set of $t$-tuples
$(s_0,s_1,\ldots,s_{t-1})$ of integers satisfying (for $0 \le j
\le t-1$) the following:
\begin{equation}\label{H}\begin{array}{l}
(1)\quad 1 \le s_j \le m,
\\
(2)\quad 0 \le ps_{j+1}-s_j \le (p-1)(m+1),
\end{array}\end{equation} with the subscripts read (mod
$t$). The set ${\mathcal H}$ was introduced in \cite{ham}, and
used in \cite{bsin} to describe the module structure of
$\Ff_q^{\Ll_1}$ under the natural action of ${\rm GL}(m+1,q)$.

For a nonconstant basis monomial
$$f(x_0,x_1,\ldots,x_m)=x_0^{b_0}\cdots x_m^{b_m},$$
in ${\mathcal M}$, we expand the exponents
$$b_i=a_{i,0}+pa_{i,1}+\cdots+p^{t-1}a_{i,t-1}\quad 0\le
a_{i,j}\le p-1$$ and let
\begin{equation}\label{lambda}
\lambda_j=a_{0,j}+\cdots+a_{m,j}.
\end{equation}
Because the total degree $\sum_{i=0}^mb_i$ is divisible by $q-1$,
there is a uniquely defined $t$-tuple $(s_0,\ldots,s_{t-1})\in
\mathcal H$ \cite{bsin} such that
$$\lambda_j=ps_{j+1}-s_j.$$ Explicitly
\begin{equation}\label{defs}
s_j=\frac1{q-1}\sum_{i=0}^m\big(\sum_{\ell=0}^{j-1}p^{\ell+t-j}a_{i,\ell}+
\sum_{\ell=j}^{t-1}p^{\ell-j}a_{i,\ell} \big)
\end{equation}
One way of interpreting the numbers $s_j$ is that the total degree
of $f^{p^i}$ is $s_{t-i}(q-1)$, when the exponent of each
coordinate $x_i$ is reduced to be no more than $q-1$ by the
substitution $x_i^q=x_i$.  We will say that $f$ is of {\it type}
$(s_0,s_1,\ldots,s_{t-1})$.

Let $c_i$ be the coefficient of $x^{i}$ in the expansion of
$(\sum_{k=0}^{p-1}x^k)^{m+1}$. Explicitly,
$$c_i=\sum_{j=0}^{\lfloor i/p\rfloor}(-1)^j{m+1 \choose j}{m+i-jp \choose m}.$$

\begin{lemma}
Let $c_i$ and $\lambda_j$ be as defined above. The number of basis
monomials in ${\mathcal M}$ of type $(s_0,s_1,\ldots,s_{t-1})$ is
$\prod_{j=0}^{t-1}c_{\lambda_j}$.
\end{lemma}

The proof of this lemma is straightforward, see \cite{csx}. For
$(s_0,s_1,\ldots ,s_{t-1})\in \mathcal H$, we will use
$c_{(s_0,s_1,\ldots ,c_{t-1})}$ to denote the number of basis
monomials in $\mathcal M$. The above lemma gives a formula for
$c_{(s_0,s_1,\ldots ,c_{t-1})}$.

\subsection{Modules of the general linear group, Hamada's formula
and the SNF of $A$}

Let $G={\rm GL}(m+1,q)$. Then $G$ acts on $\mathcal L_1$ and
$\mathcal L_d$, and $G$ is an automorphism group of the design $
(\mathcal L_1,\mathcal L_d)$. Hence each $M_i$ is an
$RG$-submodule of $R^{\mathcal L_1}$ and each $\MM_i$ is an
$\Ff_qG$-submodule of $\Ff_q^{\mathcal L_1}$. In \cite{bsin}, the
submodule lattice of $\Ff_q^{\mathcal L_1}$ is completely
determined. We will need the following result which follows easily
from the results in \cite{bsin}. To simplify the statement of the
theorem, we say that a basis monomial $x_0^{b_0}x_1^{b_1}\cdots
x_m^{b_m}$  {\it appears} in a function $f\in \Ff_q^{\mathcal
L_1}$ if when we write $f$ as the linear combination of basis
monomials, the coefficient of $x_0^{b_0}x_1^{b_1}\cdots x_m^{b_m}$
is nonzero.

\begin{theorem}\label{AB}\hspace{-0.1in}

\begin{enumerate}
\item Every $\Ff_qG$-submodule of $\Ff_q^{\Ll_1}$ has a basis
consisting of all basis monomials in the submodule.
\item Let $M$ be any $\Ff_qG$-submodule of $\Ff_q^{\mathcal L_1}$
and let $f\in \Ff_q^{\mathcal L_1}$. Then $f\in M$ if and only if
each monomial appearing in $f$ is in $M$.
\end{enumerate}
\end{theorem}

For the proof of (1), see \cite{csx}. Part (2) follows from part
(1) easily. The following is the main theorem on $\MM_1$. It was
proved by Delsarte \cite{delsarte} in 1970, and later in
\cite{glyn} and \cite{bsin}.

\begin{theorem}\label{m1}
Let $\MM_1$ be defined as above, i.e., $\MM_1$ is the dual code of
the $q$-ary (block) code of the 2-design $({\mathcal L_1},
{\mathcal L_d})$.
\begin{enumerate}
\item For any $f\in\Ff_q^{\mathcal L_1}$, we have $f\in \MM_1$ if and only if every basis monomial appearing in $f$ is in $\MM_1$.
\item Let $x_0^{b_0}x_1^{b_1}\cdots x_m^{b_m}$ be a basis monomial of type $(s_0,s_1,\ldots ,s_{t-1})$.
Then $x_0^{b_0}x_1^{b_1}\cdots x_m^{b_m}\in \MM_1$ if and only if
there exists some $j$, $0\leq j\leq t-1$, such that $s_j<d$.
\end{enumerate}
\end{theorem}

This is what Glynn and Hirschfeld \cite{glyn} called ``the main
theorem of geometric codes". As a corollary, we have

\begin{corollary}\label{hamada}\hspace{0.1in}
\begin{enumerate}
\item The dimension of $\MM_1$ is
$$\dim_{\Ff_q}\MM_1=\sum_{(s_0,s_1,\ldots ,s_{t-1})\in
\mathcal H\atop {\exists j}, s_j<d}c_{(s_0,s_1,\ldots
,s_{t-1})}.$$
\item The $p$-rank of $A$ is $${\rm rank}_p(A)=1+\sum_{(s_0,s_1,\ldots ,s_{t-1})\in
\mathcal H\atop {\forall j}, s_j\geq d}^{}c_{(s_0,s_1,\ldots
,s_{t-1})}.$$
\end{enumerate}
\end{corollary}

The rank formula in part (2) of the above corollary is the
so-called Hamada's formula.

Generalizing Theorem~\ref{m1}, we proved the following theorem in
\cite{csx}.

\begin{theorem}\label{mgeq1}
Let $\alpha\geq 1$ be an integer, and let $\MM_{\alpha}$ be
defined as above.
\begin{enumerate}
\item For any $f\in\Ff_q^{\mathcal L_1}$, we have $f\in \MM_{\alpha}$ if and only if every basis monomial appearing in $f$ is in $\MM_{\alpha}$.
\item Let $x_0^{b_0}x_1^{b_1}\cdots x_m^{b_m}$ be a basis monomial of type $(s_0,s_1,\ldots ,s_{t-1})$.
Then $x_0^{b_0}x_1^{b_1}\cdots x_m^{b_m}\in \MM_{\alpha}$ if and
only if $\sum_{j=0}^{t-1}{\rm max}\{0, d-s_j\}\geq \alpha$
\end{enumerate}
\end{theorem}

An immediate corollary is

\begin{corollary}\label{pc}
Let $0\leq \alpha\leq (d-1)t$, and let $h(\alpha,m,d+1)$ be the
multiplicity of $p^{\alpha}$ appearing as an invariant factor of
$A$. Then
$$h(\alpha,m,d+1)=\delta(0,\alpha)+\sum_{(s_0,s_1,\ldots ,s_{t-1})\in
\mathcal H \atop \sum_{j}{\rm max}\{0,
d-s_j\}=\alpha}c_{(s_0,s_1,\ldots ,s_{t-1})},$$ where
\begin{equation*} \delta(0, \alpha)= \left\{
\begin{array}{ll}
1, & \mbox{if} \; \alpha=0,\\
0, & \mbox{otherwise}. \end{array} \right.
\end{equation*}
\end{corollary}

We give some indication on how Theorem~\ref{mgeq1} was proved in
\cite{csx}. Of course Part (1) of Theorem~\ref{mgeq1} follows from
the more general result in Theorem~\ref{AB}. About Part (2) of the
theorem, if $\sum_{j=0}^{t-1}{\rm max}\{0, d-s_j\}\geq \alpha$, we
need to show that there exists a lifting of the monomial
$x_0^{b_0}x_1^{b_1}\cdots x_m^{b_m}$ to $R^{\mathcal L_1}$ that is
in $M_{\alpha}$. It turns out that the Teichm\"uller lifting
$T(x_0^{b_0}x_1^{b_1}\cdots x_m^{b_m})$ of
$x_0^{b_0}x_1^{b_1}\cdots x_m^{b_m}$ will suit our purpose. 
Indeed to show that $T(x_0^{b_0}x_1^{b_1}\cdots x_m^{b_m})\in
M_{\alpha}$, we used a theorem of Wan \cite{wan} which gives a
lower bound on the $p$-adic valuation of multiplicative character
sums. For details, we refer the reader to \cite{csx}. The other
direction of Part (2) of Theorem~\ref{mgeq1} is much more
difficult to prove. We need to prove that if $\sum_{j=0}^{t-1}{\rm
max}\{0, d-s_j\}<\alpha$, then no lifting of
$x_0^{b_0}x_1^{b_1}\cdots x_m^{b_m}$ to $R^{\mathcal L_1}$ is in
$M_{\alpha}$. We need to use the action of $G$ on $M_{\alpha}$,
Jacobi sums, and Stickelberger's theorem on Gauss sums to achieve
this. See \cite{csx} for details.

\subsection{The SNF of the 2-design in Example~\ref{affexample}}
Let ${\rm AG}(m,q)$ be the $m$-dimensional affine space over $\Ff_q$, where $q=p^t$, $p$ is a prime. Let $\cD$ be
the design in Example~\ref{affexample}, i.e., the design of the points and $d$-flats of ${\rm AG}(m,q)$. Let $A'$
be an incidence matrix of $\cD$. By viewing ${\rm AG}(m,q)$ as obtained from ${\rm PG}(m,q)$ by deleting a
hyperplane, we prove the following theorem in \cite{csx}.

\begin{theorem}\label{af}
The invariant factors of $A'$ are $p^{\alpha}$, $0\leq \alpha\leq
dt$, with multiplicity $h(\alpha,m,d+1)-h(\alpha,m-1,d+1)$, where
$h(\alpha, \cdot ,\cdot)$ is defined in Corollary~\ref{pc}.
\end{theorem}

In closing this section, we mention the following open problem.
Adopting the notation introduced at the beginning of this section,
we let $A_{d,e}$ be a (0,1)-matrix with rows indexed by elements
$Y$ of $\mathcal L_d$ and columns indexed by elements $Z$ of
$\mathcal L_e$, and with the $(Y,Z)$ entry equal to 1 if and only
if $Z\subset Y$. Note that $A_{d,1}=A$, an incidence matrix of the
2-design $(\mathcal L_1, \mathcal L_d)$. We are interested in
finding the Smith normal form of $A_{d,e}$ when $e>1$.

\begin{problem}\label{de}
Let $e>1$. What is the $p$-rank of $A_{d,e}$? And what is the SNF
of $A_{d,e}$?
\end{problem}

The first question in Problem~\ref{de} appeared in \cite{godsil},
and later in \cite{bagchi}. The $\ell$-rank of $A_{d,e}$, where
$\ell\neq p$ is a prime, is known from \cite{fy}.

\section{$p$-ranks and SNF of unitals}

A {\it unital} is a 2-$(m^3+1, m+1, 1)$ design, where $m\geq 2$.
All known unitals with parameters $(m^3+1, m+1, 1)$ have $m$ equal
to a prime power, except for one example with $m=6$ constructed by
Mathon \cite{ma}, and independently by Bagchi and Bagchi
\cite{bb}. In this section, we will only consider unitals embedded
in ${\rm PG}(2,q^2)$, i.e., unitals coming from a set of $q^3+1$
points of ${\rm PG}(2,q^2)$ which meets every line of ${\rm
PG}(2,q^2)$ in either 1 or $q+1$ points. A classical example of
such unitals is {\it the Hermitian unital} ${\mathcal U}=(\mathcal
P, \mathcal B)$, where $\mathcal P$ and $\mathcal B$ are the set
of absolute points and the set of non-absolute lines of a unitary
polarity of ${\rm PG}(2,q^2)$ respectively. Note that the order of
$\mathcal U$ is $q^2-1$. By Theorem~\ref{genrank}, only the codes
$C_p(\mathcal U)$, with $p|(q-1)$ or $p|(q+1)$, are of interest.
Furthermore, it was shown in \cite{mor} that the codes
$C_p(\mathcal U)$, with $p|(q-1)$ but $p\not{|}(q+1)$, are the
full space. So we only need to consider $C_p(\mathcal U)$ with
$p|(q+1)$. It was conjectured by Andriamanalimanana \cite{thesis}
(see also \cite{akupdate}) that ${\rm rank}_p(\mathcal
U)=(q^2-q+1)q$, if $p$ is a prime dividing $q+1$. The same
conjecture also arose in the work of Geck \cite{geck}, in which he
established a close connection between ${\rm rank}_p(\mathcal U)$
and certain decomposition numbers of the three dimensional unitary
group.

Building upon \cite{geck}, and the recent important work of
Okuyama and Waki \cite{ow} on decomposition numbers of ${\rm
SU}(3,q^2)$, Hiss \cite{hiss} determined the SNF of $\mathcal U$,
hence found the $p$-rank of $\mathcal U$ for every prime $p$.

\begin{theorem}
The invariant factors of $\mathcal U$ are
$$1^{(q^3-q^2+q)}(q+1)^{(q^2-q+1)},$$
where the exponents indicate the multiplicities of the invariant
factors. In particular, ${\rm rank}_p(\mathcal U)=q^3-q^2+q$ if
$p|(q+1)$.
\end{theorem}

The Hermitian unital is a special example of a large class of
unitals embedded in ${\rm PG}(2,q^2)$, called the Buekenhout-Metz
unitals. We refer the reader to \cite{ebert} for a survey of
results on these unitals. A subclass of the Buekenhout-Metz
unitals which received some attention can be defined as follows.

Let $q$ be an odd prime power, let $\beta$ be a primitive element
of $\Ff_{q^2}$, and for $r\in \Ff_q$ let $C_r=\{(1,y,\beta
y^2+r)\mid y\in\Ff_{q^2})\}\cup\{(0,0,1)\}$. We define
$$U_{\beta}=\cup_{r\in \Ff_q} C_r.$$
Note that each $C_r$ is a conic in ${\rm PG(2,q^2)}$, any two
distinct $C_r$ have only one point $P_{\infty}=(0,0,1)$ in common.
Hence $|U_{\beta}|=q^3+1$. It can be shown that every line of
${\rm PG}(2,q^2)$ meets $U_{\beta}$ in either 1 or $q+1$ points
(see \cite{be} and \cite{hs}). One can immediately obtain a unital
(design) ${\mathcal U}_{\beta}$ from $U_{\beta}$. We use the
points of $U_{\beta}$ as the {\it points} of $\mathcal U_{\beta}$,
and use the intersections of the secant lines with $U_{\beta}$ as
{\it blocks} to get a $2-(q^3+1,q+1,1)$ design ${\mathcal
U}_{\beta}$. Little is known about the codes of this design. As a
first step, consider the binary code $C_2({\mathcal U}_{\beta})$
of this design. The following proposition and conjecture are due
to Baker and Wantz.

To state the proposition, we use $v^{S}$ to denote the
characteristic vector of a subset $S$ in $U_{\beta}$.

\begin{proposition}[Baker and Wantz]
The vectors $v^{C_r}, r\in\Ff_q$, form a linearly independent set
of vectors in $C_2({\mathcal U_{\beta}})^{\perp}$.
\end{proposition}

\begin{proof} A binary vector $v$ lies in $C_2({\mathcal
U_{\beta}})^{\perp}$ if and only if each block of the design meets
the support of $v$ in an even number of points. If a block of
${\mathcal U_{\beta}}$ goes through $P_{\infty}$, then it meets
every $C_r$ in 2 points; if a block of ${\mathcal U}_{\beta}$ does
not go through $P_{\infty}$, then it meets every $C_r$ in either 0
or 2 points. Hence $v^{C_r}\in C_2({\mathcal U})^{\perp}$, for
every $r\in\Ff_q$. The $q$ conics $C_r$ have only one point
$P_{\infty}$ in common. Thus, $v^{C_r}$, $r\in\Ff_q$, are linearly
independent. This completes the proof.
\end{proof}

An immediate corollary is that ${\rm dim}C_2({\mathcal
U_{\beta}})^{\perp}\geq q$, hence ${\rm dim}C_2({\mathcal
U_{\beta}})\leq q^3+1-q$. Baker and Wantz made the following
conjecture.

\begin{conjecture} [Baker and Wantz]
The 2-rank  of ${\mathcal U_{\beta}}$ is $q^3+1-q$.
\end{conjecture}

Further computations done by Wantz \cite{wantz} seem to suggest
that all invariant factors of ${\mathcal U_{\beta}}$ are
$2$-powers, except for the last one, which is a 2-power times
$q+1$.

\vspace{0.2in}

\noindent{\bf Acknowledgements}: This paper was prepared while the author was visiting Caltech in January, 2004.
The author would like to thank Rick Wilson for his support.

\bibliographystyle{amsalpha}

\end{document}